\documentclass[oneside]{amsart}
\usepackage{mathrsfs}
\usepackage{graphicx} 

\newcommand{\hide}[1]{}

\usepackage{amssymb}
\usepackage{graphicx}
\usepackage{amsmath}

\usepackage{geometry}
\usepackage{caption}

\usepackage{wrapfig}

\usepackage{float}

\usepackage{bm}

\DeclareMathAlphabet{\mathdutchcal}{U}{dutchcal}{m}{n}
\SetMathAlphabet{\mathdutchcal}{bold}{U}{dutchcal}{b}{n}
\DeclareMathAlphabet{\mathdutchbcal}{U}{dutchcal}{b}{n}

\usepackage{tikz}
\usetikzlibrary{arrows, arrows.meta}

\def\textcolor#1{}

\DeclareMathOperator{\Isom}{Isom}

\newcommand{\N}{\mathbb{N}}
\newcommand{\R}{\mathbb{R}}

\renewcommand{\S}{\mathbb{S}}

\newcommand{\B}{{\bm B}}

\newcommand{\cC}{\mathcal C}

\renewcommand{\tilde}{\widetilde}
\renewcommand{\rho}{\varrho}
\renewcommand{\phi}{\varphi}

\renewcommand{\theta}{\vartheta}

\DeclareMathOperator{\Vol}{Vol}

\renewcommand{\ge}{\geqslant}
\renewcommand{\le}{\leqslant}

\usepackage[colorlinks=true,urlcolor=blue]{hyperref}

\theoremstyle{theorem}
\newtheorem{theorem}{Theorem}[section]
\newtheorem{lemma}[theorem]{Lemma}

\newtheorem{MainTheorem}{Theorem}

\newtheorem*{theoremAF}{The Alexandrov--Fenchel Inequality}

\newtheoremstyle{Intro}
{}
{}
{\itshape}
{}
{\scshape}
{.}
{.5em}
{}

\theoremstyle{Intro}

\theoremstyle{definition}
\newtheorem{definition}[theorem]{Definition}

\theoremstyle{remark}
\newtheorem*{remark}{\textsc{Remark}}

\newcounter{reminder}

\newtheoremstyle{claim}
  {}
  {}
  {\itshape}
  {0pt}
  {\scshape}
  {.}
  { }
  {\thmname{#1}\thmnumber{ #2}\thmnote{ (#3)}}

\theoremstyle{claim}

\numberwithin{equation}{section}

\title{A solution to Bezdek's conjecture}

\author[Kostiantyn Drach]{Kostiantyn Drach}
\author[Kateryna Tatarko]{Kateryna Tatarko}
\thanks{ The authors would like to thank Maud Szusterman and Dmitry Ryabogin for valuable discussions. The first author is partially supported by Departament de Recerca i Universitats de la Generalitat de Catalunya (2021 SGR 00697) and by Agencia Estatal de Investigaci\'on Grant PID2023-147252NB-I00 and through the Severo Ochoa and Mar\'ia de Maeztu Program for Centers and Units of Excellence in R\&D (CEX2020-001084-M). The second author is partially supported by NSERC Discovery Grant number 2022-02961. }

\date{}

\address{Universitat de Barcelona, Gran Via de les Corts Catalanes, 585, 08007 Barcelona, Spain}
\address{Centre de Recerca Matem\`atica, Edifici C, Carrer de l'Albareda, 08193 Bellaterra, Barcelona, Spain}
\email{kostiantyn.drach@ub.edu}
\address{Department of Pure Mathematics, University of Waterloo, Waterloo, ON, N2L 3G1, Canada}
\email{ktatarko@uwaterloo.ca}

\begin{document}

\begin{abstract}
For a given $\lambda >0$, a convex body in $\R^n$ is \emph{$\lambda$-convex} if it is the intersection of (finitely or infinitely many) balls of radius $1/\lambda$. In this note, we show that among all $\lambda$-convex bodies in $\R^n$, $n \ge 2$, with a given inradius, the $\lambda$-convex lens (i.e., the intersection of two balls of radius $1/\lambda$) has the largest mean width. This gives an affirmative answer to the conjecture of K.~Bezdek. Under an additional symmetry assumption on $\lambda$-convex bodies, we resolve the analogous inradius conjecture of Bezdek for arbitrary intrinsic volumes.

We also establish an answer to the corresponding conjecture of K.~Bezdek about the circumradius. In particular, we prove that the $\lambda$-convex spindle (i.e., the intersection of all balls of radius $1/\lambda$ containing a given pair of points) is the unique minimizer of the mean width among all $\lambda$-convex bodies with a fixed circumradius.
\end{abstract}

\maketitle

\section{Introduction}

A convex body in $\R^n$ is a compact convex subset with non-empty interior. For a given $\lambda >0$, a convex body $K \subset \R^n$ is called \emph{$\lambda$-convex} if for every point $p \in \partial K$ there exists a neighborhood $U_p$ and a ball $B_{\lambda, p}$ of radius $1/\lambda$ such that $p \in \partial B_{\lambda,p}$ and $U_p \cap K\subset B_{\lambda,p}$ (see Figure~\ref{Fig:Def}). By Blaschke's rolling theorem \cite{Bla56} (see also \cite{DrBla} and references therein), every $\lambda$-convex body can be represented as an intersection of (finitely or infinitely many) balls of radius $1/\lambda$. A \emph{$\lambda$-convex lens} in $\R^n$ is the intersection of two balls of radius $1/\lambda$.

The class of $\lambda$-convex bodies appears naturally in connection with the Kneser--Poulsen conjecture \cite{Bezdek2008}, the study of the illumination conjecture \cite{Bezdek2012}, approximation of convex bodies \cite{FV, NV}; see also \cite{AAF, AACF, SWY} for further developments on the geometry of $\lambda$-convex bodies and related questions. Recently, the study of \emph{reverse} isoperimetric problems in the class of $\lambda$-convex bodies garnered a lot of attention. In particular, the analogue of Ball's reverse isoperimetric inequality for $\lambda$-convex bodies was obtained in \cite{FKV, BD} in $\mathbb{R}^2$ and in \cite{DT} in $\mathbb{R}^3$.   We refer the reader to the nice surveys \cite{BLN, BLNP} about topics in $\lambda$-convexity.

The classical Steiner formula asserts that for every $\varepsilon >0$, 
\begin{equation}\label{Eq:SteinerForm}
\Vol(K+ \varepsilon \B) = \sum\limits_{j = 0}^n \kappa_{n-j} V_j(K) \varepsilon^{n-j}
\end{equation}
where $\Vol$ is the Lebesgue measure on $\mathbb{R}^n$, $\B$ is the unit Euclidean ball in $\mathbb{R}^n$, $\kappa_{n-j}$ is the volume of the unit Euclidean ball in $\mathbb{R}^{n-j}$, and $+$ denotes the Minkowski addition. The coefficients $V_j(K)$ are called the \emph{intrinsic volumes} of $K$. It is known that $V_n(K)$ is the volume of $K$, $V_{n-1}(K)$ is (up to the scalar) the surface area of $K$, $V_{1}(K)$ is (up to a scalar) the mean width of $K$, and $V_{0}(K) = \kappa_{n}$ is the volume of $\B$. We refer readers to \cite{Sch} for a comprehensive discussion about intrinsic volumes.

The \emph{inradius} $r(K)$ of a convex body $K$ is the radius of
a largest ball contained in $K$. In this paper, we answer the conjecture posed by K. Bezdek in \cite[Conjecture 5]{BezdekConj} about minimizing the inradius in the class of $\lambda$-convex bodies with fixed mean width.
We remark that this result supplements the analogous statements for the volume due to Bezdek \cite{BezdekConj}, and for the surface area due to the authors \cite{DT}.

\begin{MainTheorem}[Reverse inradius inequality for mean width]
\label{Thm:Main}
Let $n\ge 2$ and let $\lambda > 0$. Let $K \subset \R^n$  be a $\lambda$-convex body and $L \subset \R^n$ be a $\lambda$-convex lens. If $r(K)=r(L)$,
then
\begin{equation}
\label{Eq:Goal}
V_1(K) \le V_1(L),
\end{equation}
with equality if and only if $K$ is a $\lambda$-convex lens.
\end{MainTheorem}

For $n=2$, this result was obtained by Milka in \cite{MilInradius}, and it was later generalized to all $2$-dimensional Alexandrov spaces of curvature bounded below by the first author in \cite{Dr}.

The circumradius $R(K)$ of a convex body $K$ is the radius of the smallest ball that contains $K$. K. Bezdek \cite[Conjecture 10]{BezdekConj} conjectured that a \emph{$\lambda$-convex spindle} has the largest circumradius among all $\lambda$-convex bodies in $\mathbb{R}^n$ with fixed mean width. A \emph{$\lambda$-convex spindle} is the intersection of all balls that contain a pair of given points. These points, called \emph{the vertices} of the spindle, must be necessarily at a distance at most $2/\lambda$ apart. Spindles appear as solutions to various optimization problems for $\lambda$-convex domains (see e.g., \cite{BezdekConj, BorDr13, BM, DrSpindle}).

As an application of Theorem~\ref{Thm:Main} and the use of the $\lambda$-duality in the class of $\lambda$-convex bodies, we obtain a solution to the partial case of Bezdek's conjecture \cite[Conjecture 10]{BezdekConj} for the mean width.

\begin{MainTheorem}[Reverse circumradius inequality for mean width]
\label{Thm:Main2}
Let $n\ge 2$ and let $\lambda > 0$. Let $K \subset \R^n$ be a $\lambda$-convex body and $S \subset \R^n$ be a $\lambda$-convex spindle. If $R(K)=R(S)$,
then
\begin{equation}
\label{Eq:Meanwidth}
V_1(K) \ge V_1(S), 
\end{equation}
with equality if and only if $K$ is a $\lambda$-convex spindle.
\end{MainTheorem}

We note that the inequality of Linhart \cite{Li1977} gives a lower bound for $V_1(K)$ in terms of the circumradius  of $K$ for any convex body $K \subset \mathbb{R}^n$, that is,
$$
V_1(K) \ge 2 R(K)
$$
with equality if and only if $K$ is a segment (see also \cite[Theorem 3.3]{BH}). In view of remark in \cite[Remark 11]{BezdekConj}, Theorem~\ref{Thm:Main2} extends Linhart's inequality  from the classical convexity to $\lambda$-convexity.

If a $\lambda$-convex body $K \subset \mathbb{R}^n$ has \emph{inball symmetries} (that is, $K$ possesses an additional symmetry assumption defined in Section~\ref{ThmCProof}), we resolve K. Bezdek's conjecture \cite[Conjecture 5]{BezdekConj} for any $j$-th intrinsic volume.

\begin{MainTheorem}[Reverse inradius inequality for intrinsic volumes under symmetry assumption]
\label{Thm:MainSym}
Let $n\ge 2$ and let $\lambda > 0$. Let $K \in \R^n$ be a $\lambda$-convex body that has inball symmetries and $L \subset \R^n$ be a $\lambda$-convex lens. If $r(K)=r(L)$,
then
\begin{equation}
\label{Eq:Goal}
V_j(K) \le V_j(L) \qquad \text{for any} \quad j \in \{1,\ldots, n\},
\end{equation}
with equality if and only if $K$ is a $\lambda$-convex lens.
\end{MainTheorem}

The paper is organized as follows. In Section~\ref{background}, we introduce basic notions and provide the needed background. In Sections~\ref{ThmAProof} and \ref{ThmBProof}, we prove Theorem~\ref{Thm:Main} and Theorem~\ref{Thm:Main2}, respectively. Finally, imposing additional symmetry constraints on $\lambda$-convex bodies, we prove Theorem~\ref{Thm:MainSym} in Section~\ref{ThmCProof}.

\section{Preliminaries}\label{background}

In this section, we aim to recall some general background from convex geometry and $\lambda$-convexity.

\subsection{Mixed volumes.} We refer the readers to \cite{Sch} for a detailed study of mixed volumes and basic notions from convexity.

We work in the Euclidean space $\R^n$ with scalar product $\langle \cdot, \cdot\rangle.$ We denote the unit ball centered at the origin by $\B$, and its boundary unit sphere as $\S^{n-1} = \partial \B$.

Let $K_1, \dots, K_m$ be convex bodies in $\R^n$. The fundamental theorem of Minkowski describes the behavior of the volume of the Minkowski sum of convex bodies:
\begin{equation}\label{Eq:MinkForm}
    \Vol(t_1 K_1 + \dots + t_m K_m) =  \sum\limits_{i_1, \dots, i_n = 1}^m V(K_{i_1}, \dots, K_{i_n}) t_{i_1} \dots t_{i_n}
\end{equation}
for any $t_1, \dots, t_m \geq 0.$ The coefficients $V(K_{i_1}, \dots, K_{i_n})$  in this polynomial are called {\it mixed volumes} which are important geometric functionals associated with convex bodies. These coefficients are nonnegative, continuous, and multilinear in their arguments. Moreover, mixed volumes are invariant under permutations, that is, if $\tau$ is a permutation on $\{1, \ldots, n\}$, then
\begin{equation}
    \label{Eq:Sym}
    V(K_1, \ldots, K_n) = V(K_{\tau(1)}, \ldots, K_{\tau(n)}).
\end{equation}

As a special case, \eqref{Eq:MinkForm} implies the Steiner's formula \eqref{Eq:SteinerForm}. The coefficients $V(K, \dots, K, \B, \dots, \B)$ are called quermassintegrals $W_j(K)$ of a convex body $K$ and defined as   
\[
W_j(K) = V(\underbrace{K, \ldots, K}_{n-j}, \underbrace{\B, \ldots, \B}_{j})
\]
for $j = 0, 1, \dots, n.$
We will use their weighted versions $V_j(K)$ which are called {\it intrinsic volumes} of~$K$:
\begin{equation}\label{Eq:IntrVol}
V_j(K) = \frac{{n \choose j}}{\kappa_{n-j}} W_{n-j}(K) = \frac{{n \choose j}}{\kappa_{n-j}} V(\underbrace{K, \ldots, K}_{j}, \underbrace{\B, \ldots, \B}_{n-j}).
\end{equation}
It is known that $V_n(K)$ is the $n$-dimensional volume of $K$, $2 V_{n-1}(K)$ is the surface area  of $K$, $\frac{2 \kappa_{n-1}}{n \kappa_n}V_{1}(K)$ is the mean width of $K$, and finally, $V_0(K) = \kappa_n$ is just the volume of the unit Euclidean ball $\B$ in $\R^n$.

The support function of a convex body $K$ is defined as $h_K(x) = \max\limits_{y \in K} \langle x, y \rangle$ for every $x \in \R^n$. The mixed area measure $S(K_2, \ldots, K_n)(\cdot)$ of $K_2, \dots, K_n$ is a unique non-negative measure on $\S^{n-1}$ such that
\[
V(K_1, \ldots, K_n) = \frac1n \int_{\mathbb S^{n-1}} h_{K_1}(u) \, dS(K_2, \ldots, K_n)(u).
\]

Note that $S(\B, \ldots, \B)$ is the standard measure $du$ on $\S^{n-1}$, and more generally, $S(K, \ldots, K)$ is the surface area measure on $\partial K$.

The mixed area measure $S(K_2, \ldots, K_n)(\cdot)$ is invariant under ambient isometries applied simultaneously to all convex bodies:
\begin{equation}
    \label{Eq:Isometry}
    S(K_2, \ldots, K_n)(\cdot) = S(\sigma K_2, \ldots, \sigma K_n)(\cdot), \quad \text{for any orientation-preserving isometry } \sigma \colon \R^n \to \R^n.
\end{equation}
Furthermore, if each of the bodies $K_2, \ldots, K_n$ is rotationally symmetric with respect to an axis parallel to a given unit direction $\ell$, then the measure $S(K_2, \ldots, K_n)(\cdot)$ is rotationally invariant in the following sense. Let $p \in \S^{n-1}$ be the point that corresponds to $\ell$, and let $A \subset \S^{n-1}$ be a set measurable with respect to $S(K_2, \ldots, K_n)(\cdot)$. Then for every rotation $\sigma$ on $\S^{n-1}$ around $p$, we have 
\begin{equation}
\label{Eq:Rotation}
S(K_2, \ldots, K_n)(A) = S(K_2, \ldots, K_n)(\sigma A).
\end{equation}

The \emph{width} of a convex body $K \subset \R^n$ in the direction $u \in\S^{n-1}$ is defined as $h_K(u) + h_K(-u)$. The \emph{mean width} $w(K)$ of $K$ is the average of the widths of $K$, that is,
\begin{equation}
\label{Eq:MeanWidth}
w(K) = \frac{1}{|\S^{n-1}|}\int_{\S^{n-1}} \left(h_K(u) + h_K(-u)\right) \,du = \frac{2}{|\S^{n-1}|} \int_{\S^{n-1}} h_K(u) \,du,
\end{equation}
where $|\S^{n-1}|$ is the $(n-1)$-dimensional measure of $\S^{n-1}$. Thus, we can write $V_1(K)$ as
\begin{equation}
\label{Eq:MeanWidth}
V_1(K) =  \frac{1}{\kappa_{n-1}}  \int_{\S^{n-1}} h_K(u) \,du = \frac{ n \kappa_n}{2 \kappa_{n-1}} w(K),
\end{equation}
where we used that $|\S^{n-1}| = n \kappa_n$.

The most fundamental inequality that relates mixed volumes is the classical Alexandrov--Fenchel inequality. 

\begin{theoremAF}
    Let $K_1, K_2, \ldots, K_n \subset \R^n$ be convex bodies. Then
    \[
    V(K_1, K_2, K_3, \dots, K_n)^2 \ge V(K_1, K_1, K_3, \dots, K_n) \cdot V(K_2, K_2, K_3, \dots, K_n).
    \]
\end{theoremAF}


\subsection{$\lambda$-convexity.} In this subsection, we collect some definitions and properties related to $\lambda$-convex bodies.

Let $\lambda>0$. A {\it $\lambda$-convex polytope} is a $\lambda$-convex body given as the intersection of finitely many balls of radius $1\slash\lambda$. An $(n-1)$-dimensional {\it facet} of a  $\lambda$-convex polytope  $P$ is the intersection of the corresponding ball of radius $1\slash\lambda$  with the boundary of $P$. Two facets are either disjoint or intersect. In the latter case, their intersection is called a (lower-dimensional) face. The faces
of dimension $0$ and $1$ are called vertices and edges of the polytope $P$, respectively. The $\lambda$-convex lens is an example of a $\lambda$-convex polytope.

\begin{figure}[h!] 
    \centering
     \includegraphics[scale=0.3]{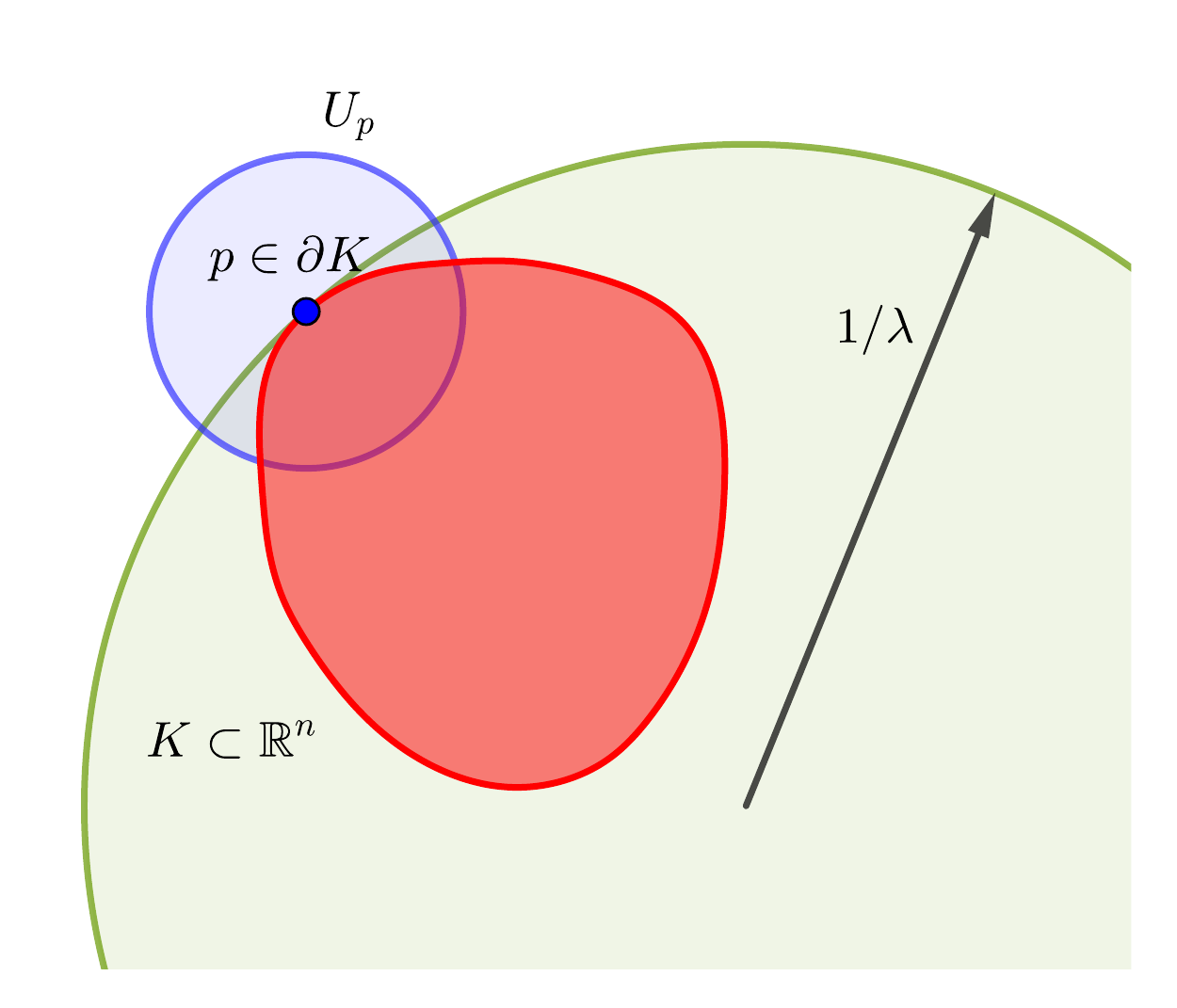}
   \caption{The definition of a $\lambda$-convex body}
    \label{Fig:Def}
\end{figure}

We also recall that a convex body $K \subset \R^n$ is called \emph{$\lambda$-convex} if for every point $p \in \partial K$ there exists a neighborhood $U_p$ and a ball $B_{\lambda, p}$ of radius $1/\lambda$ such that $p \in \partial B_{\lambda,p}$ and $U_p \cap K\subset B_{\lambda,p}$ (see Figure~\ref{Fig:Def}). When $\partial K$ is at least $C^2$-smooth, this condition is equivalent to requiring that the principal curvatures $k_i(p)$, $i\in\{1, \ldots, n-1\}$, with respect to the inward-pointing normal at every $p\in \partial K$, satisfy $k_i(p) \ge 1$. 

For every $\lambda$-convex body $K$, we define the {\it $\lambda$-convex dual} $K^\lambda$ as
$$
K^\lambda = \underset{x \in K}{\bigcap} \left(\frac1{\lambda} \B + x\right).
$$
It follows from this definition that $K^\lambda$ is also a $\lambda$-convex body. Moreover, $\lambda$-duality is an order-reversing involution, that is, $(K^\lambda)^\lambda = K$ and $K^\lambda \subseteq M^\lambda$ if $M \subseteq K$ for any $\lambda$-convex bodies $K$ and $M$. The $\lambda$-duality has been recently investigated in \cite{AACF, BLN, FKV}.

It is known (see, for example, \cite[Proposition 2.2]{FKV}) that for a $\lambda$-convex body $K$ the following relation holds
\[
K + (-K^\lambda)=\frac1{\lambda}\B.
\]
Equivalently, 
\begin{equation}
\label{Eq:Support}
h_K(u)+h_{K^\lambda}(-u)= \frac1{\lambda}, \qquad \text{for any } u\in \S^{n-1}.
\end{equation}
The relation~\eqref{Eq:Support} immediately implies that 
\begin{equation}
\label{Eq:Inrad}
r(K) + R(K^\lambda) = \frac1{\lambda}.   
\end{equation}
Finally, if $L$ is a $\lambda$-convex lens, then $L^\lambda = S$ is a $\lambda$-convex spindle. We refer the reader to \cite[Section 3.1]{BLN} for the detailed discussion of the $\lambda$-duality.

\begin{remark}
    We note that  $\lambda$-duality is different from the spherical duality between $\lambda$-convex and $\lambda$-concave domains in $\S^2$, which was used in \cite{BorDr15_1} to relate the reverse isoperimetric problems in these two classes of bodies. 
\end{remark}

\section{Proof of Theorem \ref{Thm:Main}}\label{ThmAProof}

Without loss of generality, we can assume that $\lambda=1$. Let $B$ be an inscribed ball for $K$. We can assume that this ball is centered at the origin and has the radius $r(K) = r(L) \le 1$.

Since $B$ is the inscribed ball for $K$, the intersection $\partial K \cap \partial B$ contains a finite set of points $T = \{p_0, p_1,\ldots\}$, with $|T| \le n+1$, that does not lie in any of the open hemispheres of $\partial B$ (this statement is folklore, see e.g., \cite[Proposition 2.9]{DT}; it can also be deduced from the Carath\'eodory theorem).

Let $\tilde K$ be a $1$-convex polytope given by the intersection of the unit balls tangent to $B$ at the points in $T$ and containing $B$ in their closures. By Blaschke's Rolling Theorem, $K \subseteq \tilde K$. Furthermore, by \cite[Proposition 2.9]{DT}, $B$ is the inscribed ball for $\tilde K$. Thus, 
\[
    r(\tilde K) = r(K) \quad \text{and} \quad K \subseteq \tilde K.
\]
Hence, 
\[
  V_1(K) \le V_1(\tilde K)
\]
by the monotonicity of the intrinsic volumes, with equality if and only if $\tilde K = K$.

Therefore, to establish Theorem \ref{Thm:Main}, it is enough to consider the case when $K$ is a $1$-convex polytope with the property that all of its facets touch the inscribed ball at the points in the set $T = \{p_0, p_1, \ldots,p_{|T|}\}$, $|T| \le n+1$. 

For the rest of this section, we further assume that the lens $L$ is a $1$-convex lens centered at the origin (and hence, sharing the same inscribed ball with $K$) and rotated so that it touches the inscribed ball in one of the points in $T$, say, $p_0$.  

Let $\mathcal C := (\underbrace{\B, \ldots, \B}_{n-1})$ be an $(n-1)$-tuple of unit Euclidean balls $\B$. Then our goal is to show 
\begin{equation}
\label{Eq:Go0}
V(K, \mathcal{C}) = \frac{1}{n}\int \limits_{\mathbb S^{n-1}} h_K(u) du \le \frac{1}{n}\int\limits_{\mathbb S^{n-1}} h_L(u) du = V(L, \mathcal{C}),
\end{equation}
with equality if and only if $K$ is a lens.

Let $F_i$ be a facet of the $1$-convex polytope $K$. Denote by $p_i$ the point of tangency of $F_i$ and the inscribed ball $B$. Let $\mathcal N_{F_i} \subset  \mathbb S^{n-1}$ be the radial projection of $F_i$ onto the unit sphere $\mathbb S^{n-1}$ with respect to the origin. Equivalently, if $Q_1, \ldots, Q_{|T|} \subset \partial B$ are the Voronoi cells of points $p_1, \ldots, p_{|T|} \in \partial B$, then each $\mathcal N_{F_i}$ is the radial projection of the corresponding $Q_i$ onto the unit sphere $\S^{n-1}$.

Assume that the $1$-convex lens $L$ is placed so that its center is at the origin and one of its facets, which we call $F$, touches the inscribed ball $B$ at the point $p_i$. If $\mathcal N_F \subset \mathbb S^{n-1}$ is a similar radial projection of $F$, then we have
    \begin{equation}
    \label{Eq:Inclusion}
    F_i \subseteq F, \quad \mathcal N_{F_i} \subseteq \mathcal N_{F}, \quad \text{with simultaneous equality if and only if }\quad F_i = F.
    \end{equation}
    The first inclusion follows from Blaschke's rolling theorem \cite{Bla56}, and the second inclusion follows from the symmetry and the fact that $F_i$ and $F$ are regions on the unit sphere.

    \begin{lemma}
    \label{Claim:M}
        \begin{equation}
        \label{Eq:M0}
        \frac{\int \limits_{\mathcal N_{F_i}} h_K(u) \, du}{\int \limits_{\mathcal N_{F_i}} du} \le \frac{\int\limits_{\mathcal N_F} h_L(u) \, du}{\int \limits_{\mathcal N_F} du} = \frac{V(L, \mathcal{C})}{V(\B,\mathcal{C})}.
    \end{equation}
    Moreover, the equality holds if and only if $\mathcal N_F = \mathcal N_{F_i}$ (and thus, if and only if $F = F_i$).
    \end{lemma}

Now, assume that Lemma~\ref{Claim:M} holds for every $i \in \{1, \ldots, |T|\}$. Then, summing over all facets of $K$, we obtain
\begin{equation*}
    \begin{aligned}
        V(K, \cC) &= \frac{1}{n}\int \limits_{\mathbb S^{n-1}} h_K(u) \, du = \frac{1}{n}\sum_i \int \limits_{\mathcal N_{F_i}} h_K(u) \, du  \\
        &\le \frac{V(L, \cC)}{V(\B,\mathcal{C})} \cdot \frac{1}{n} \sum_i \int \limits_{\mathcal N_{F_i}} du = \frac{V(L, \cC)}{V(\B,\mathcal{C})} \cdot \frac{1}{n}\int \limits_{\mathbb S^{n-1}} du = V(L, \cC), 
    \end{aligned}
\end{equation*}
which yields~\eqref{Eq:Go0} and finishes the proof of Theorem~\ref{Thm:Main}. In the next subsection, we prove Lemma~\ref{Claim:M}. 

\subsection{Proof of Lemma~\ref{Claim:M}} First observe that by~\eqref{Eq:Inclusion},
\[
h_K(u) \le h_L(u), \quad \forall u \in \mathcal N_{F_i}.
\]
Therefore, to prove \eqref{Eq:M0}, it is enough to establish
\begin{equation}
    \label{Eq:M2}
    \frac{\int \limits_{\mathcal N_{F_i}} h_L(u) \, du}{\int \limits_{\mathcal N_{F_i}} du} \le \frac{\int\limits_{\mathcal N_F} h_L(u) \, du}{\int \limits_{\mathcal N_F} du}.
\end{equation}

Choose the coordinates $(t, \theta)$ in $\S^{n-1}$ centred at $\frac1{r(K)} p_i$ such that a point $u \in \S^{n-1}$ can be written as
\[
u = \cos (t)\frac{p_i}{r(K)} + \sin (t) \, \theta,
\]
where $\theta \in \S^{n-2}$ (more precisely, $\theta \in p_i^\perp \cap\S^{n-1}$) and $t$ is the angle between $u$ and direction $op_i$. Then the surface area measure on $\S^{n-1}$ is $du = \left(\sin (t)\right)^{n-2} dt \, d\theta$,  and we set $q(t) := \left(\sin (t)\right)^{n-2}$ for $t \in [0,\pi/2]$. The function $q$ is strictly monotone increasing on its domain of definition. Also, denote the spherical radial function of $\mathcal N_{F_i}$ in the direction $\theta$ by 
$$
\varphi (\theta) := \max\limits_{t \in [0,{\pi}/2]}\left\{ \cos(t) \frac{p_i}{r(K)} + \sin (t) \, \theta \in \mathcal N_{F_i} \right\}.
$$
By rotational symmetry, the support function $h_L(u) = h_L(t)$ is independent of $\theta$. Furthermore, the function $t \mapsto h_L(t)$ is strictly monotone increasing on its domain of definition $t \in [0, \pi/2]$. 

Then we can rewrite
$$
\int\limits_{\mathcal N_{F_i}} h_L(u) \, du = \int\limits_{\S^{n-2}} \int\limits_0^{\varphi(\theta)} h_L(t) q(t) dt \,d\theta \qquad \text{and} \qquad \int \limits_{\mathcal N_{F_i}} du =  \int\limits_{\S^{n-2}} \int\limits_0^{\varphi(\theta)} q(t) dt \,d\theta.
$$
Also, let $\mathcal{F}: = {\int_0^{{\pi}/2} h_L(t) q(t) dt}\,/\,{\int_0^{{\pi}/2} q(t) dt}$ and consider 
\begin{equation}
\mathcal R[\varphi(\theta)] := \int\limits_0^{\varphi(\theta)} \left(h_L(t) - \mathcal{F} \right) q(t) \,dt,
\end{equation}
where $\varphi(\theta) \in [0, {\pi}/2]$. Note that $\mathcal R[{\pi}/2] =0$.

Since $t \mapsto \big(h_L(t) - \mathcal{F} \big) q(t)$  is a continuous and strictly increasing function on $[0, {\pi}/2]$, we obtain that $\mathcal R[\varphi(\theta)] \le 0$ for all possible $\varphi(\theta) \in [0, {\pi}/2]$ with equality if and only if $\varphi(\theta) = {\pi}/2$ (see \cite[Lemma 4.5]{DT}). Hence,
\begin{align*}
\int\limits_{\mathcal N_{F_i}} h_L(u) \, du &= \int\limits_{\S^{n-2}} \int\limits_0^{\varphi(\theta)} h_L(t) q(t) dt \,d\theta 
\le  \mathcal{F}  \int\limits_{\S^{n-2}}   \int\limits_0^{\varphi(\theta)}   q(t) dt \,d\theta \\
&= \frac{\int\limits_{\S^{n-2}} \int\limits_0^{\frac{\pi}2} h_L(t) q(t) dt \,d\theta }{\int\limits_{\S^{n-2}} \int\limits_0^{\frac{\pi}2} q(t) dt \,d\theta }  \int\limits_{\S^{n-2}} \int\limits_0^{\varphi(\theta)}   q(t) dt \,d\theta =  \frac{\int\limits_{\mathcal N_F} h_L(u) \, du}{\int \limits_{\mathcal N_F} du}  \int \limits_{\mathcal N_{F_i}} du,
\end{align*}
as required.

\section{Proof of Theorem~\ref{Thm:Main2}}\label{ThmBProof}

Let $K$ be a $\lambda$-convex body with circumradius $R(K)$, and let $S$ be a $\lambda$-convex spindle such that $R(S) = R(K)$.  Consider the $\lambda$-dual body $K^\lambda$ and the $\lambda$-convex lens $L := S^\lambda$.

By \eqref{Eq:MeanWidth} and  \eqref{Eq:Support},
\begin{equation}
\label{Eq:MW}
V_1(K) + V_1(K^\lambda) = \frac1\lambda V_1(\B).
\end{equation}
Also, $r(K^\lambda) = \frac1{\lambda} - R(K) = \frac1{\lambda} -R(S)=r(S^\lambda) = r(L)$ by \eqref{Eq:Inrad}. 

Theorem~\ref{Thm:Main} (applied for $j=1$) implies that 
\begin{equation} \label{Eq:ProofBAppl}
V_1(K^\lambda) \le V_1(L)
\end{equation}
with equality if and only if $K^\lambda$ is a $\lambda$-convex lens $L$. Thus, \eqref{Eq:MW} and \eqref{Eq:ProofBAppl} yield that 
\[
V_1(K) = \frac1\lambda V_1(\B)-V_1(K^\lambda) \ge \frac1\lambda V_1(\B) - V_1(L) = V_1(L^\lambda) = V_1(S)
\]
with equality if and only if $K$ is a $\lambda$-convex spindle $S$.
This completes the proof of Theorem~\ref{Thm:Main2} (including the equality case).

\section{Proof of Theorem \ref{Thm:MainSym}} \label{ThmCProof}

For intrinsic volumes different from the volume, the surface area, and the mean width, our method allows us to derive a result similar to our Theorem~\ref{Thm:Main} under the additional symmetry assumption as follows. Let $i(K) \in \partial K$ be the set of points at which the inscribed ball of $K$ touches the boundary $\partial K$. This set can be finite or infinite. For example, if $L$ is a lens, then $i(L)$ consists of two points (symmetric with respect to the center of $L$). If $S$  is a spindle, then $i(S)$ is a circle. 

Let $\Isom^+(\R^n)$ be the group of orientation-preserving isometries in $\R^n$ and  $\mathscr{O}(K) < \Isom^+(\R^n)$ be the maximal subgroup of rotations around the center $o$ of the inscribed ball to $K$ that fixes $i(K)$. In other words, a rotation $\sigma \in \Isom^+(\R^n)$ about $o$ lies in $\mathscr{O}(K)$ if 
\[
i(K) = \sigma\left(i(K)\right) := \big\{\sigma(x) \,|\, x \in i(K)\big\}.
\]

\begin{definition}[Inball symmetries]
\label{Def:Sym}
We say that $K$ has \emph{inball symmetries} if $\mathscr{O}(K)$ is not equal to the identity group and it acts transitively on the set $i(K)$: for every $x, y \in i(K)$ there exists $\sigma \in \mathscr{O}(K)$ such that $\sigma(x) = y$.
\end{definition}

Note that every lens has inball symmetries, and furthermore, $L$ is invariant under the action of $\mathscr{O}(L)$. However, a general convex body $K$ need not be invariant under the action of $\mathscr{O}(K)$. 

\subsection{Proof of Theorem~\ref{Thm:MainSym}}
The proof begins in a similar way to the proof of Theorem~\ref{Thm:Main}. Starting with $K$, we can construct a $1$-convex polytope $\tilde K$ such that all of its faces are tangent to the inscribed ball $B$ along the finite set $T$, and $V_j(K) \le V_j(\tilde K)$ (with equality if and only if $\tilde K = K$) using the monotonicity of the intrinsic volumes. We assume that $B$ is centered at the origin. 

Since $K$ has inball symmetries, we can further choose $T$ (if there is a choice) so that $\tilde K$ has inball symmetries as well. Note that $\mathscr{O}(\tilde K)$ is a subgroup of $\mathscr{O}(K)$, not necessarily proper. Since $\tilde K$ is a $1$-convex polytope that is completely determined by $i(\tilde K)$, we obtain that 
\begin{equation}
\label{Eq:Symmetry}
    \sigma(\tilde K) = \tilde K, \quad \quad \forall \sigma \in \mathscr{O}(\tilde K),
\end{equation}
i.e., not only are the touching points invariant under the action of $\mathscr{O}(\tilde K)$, but the polytope itself. 

Therefore, as in the proof of Theorem~\ref{Thm:Main}, in order to establish Theorem~\ref{Thm:MainSym}, we can assume that $K$ is a $1$-convex polytope that satisfies the \textbf{Standing Assumptions}: the faces of $K$ touch the inscribed ball in some set $T = \{p_0, p_1, \ldots,p_{|T|}\}$, $|T| \le n+1$, and $K$ satisfies~\eqref{Eq:Symmetry} with $\tilde K = K$. 

Also, let $L$ be a $1$-convex lens centered at the origin (and hence, sharing the same inscribed ball with $K$). 

\medskip

For a convex body $M$ and $t \in \{1,\ldots, n\}$, we will use a calligraphic version of $M$, i.e., $\mathscr M_t$, to denote a $t$-tuple of copies of $M$, i.e., $\mathscr M_t = (M, \ldots, M)$. In this way, for example, $V(\mathscr M_t,\mathscr C_{n-t}) =  V(\underbrace{M, \ldots, M}_{t}, \underbrace{C, \ldots, C}_{n-t})$. We extend this for $t=0$ by setting $V(\mathscr M_0,\mathscr C_{n}) =  V( \underbrace{C, \ldots, C}_{n})$.

\begin{lemma}
\label{Claim:1}
Let $K$ be the $1$-convex polytope that satisfies the Standing Assumption, $L$ be the $1$-convex lens, and $\B$ be the unit Euclidean ball in $\R^n$ as above. Then for every $s \in \{1, \ldots, n\}$ and $t \in \{0, \ldots, s-1\}$, 
\begin{equation}
    \label{Eq:Step1}
    \frac{V(\mathscr K_{s-t}, \mathscr L_{t}, \mathscr B_{n-s})}{V(\mathscr K_{s-t-1}, \mathscr L_{t+1}, \mathscr B_{n-s})} \le 1.
\end{equation}

Moreover, the equality holds if and only if $K$ is a lens.
\end{lemma}

This lemma immediately yields Theorem~\ref{Thm:MainSym}. Indeed, to show $V_j(K) \le V_j(L)$ for a given $j \in \{1, \ldots, n\}$, in view of \eqref{Eq:IntrVol}, we need to show that
\[
\frac{V(\mathscr K_{j}, \mathscr B_{n-j})}{V(\mathscr L_{j}, \mathscr B_{n-j})} \le 1.
\]
The following chain of the inequalities is obtained by applying Lemma~\ref{Claim:1} for $s=j$ and $t\in\{0, \ldots, j-1\}$
\[
\frac{1}{V(\mathscr L_{j}, \mathscr B_{n-j})} \overset{t=j-1}{\le} \frac{1}{V(\mathscr K_1, \mathscr L_{j-1}, \mathscr B_{n-j})}  \overset{t=j-2}{\le} \ldots \overset{t=1}{\le} \frac{1}{V(\mathscr K_{j-1}, \mathscr L_{1},  \mathscr B_{n-j})} \overset{t=0}{\le} \frac{1}{V(\mathscr K_{j}, \mathscr B_{n-j})}.
\]
Multiplying the above by $V(\mathscr K_{j}, \mathscr B_{n-j})$, we get
\[
\frac{V(\mathscr K_{j}, \mathscr B_{n-j})}{V(\mathscr L_{j}, \mathscr B_{n-j})} \le 1
\]
as needed. The equality case follows from the equality case in Lemma~\ref{Claim:1}. This finishes the proof of Theorem~\ref{Thm:MainSym}. The remainder of the section will be devoted to the proof of Lemma~\ref{Claim:1}.

\subsection{Proof of Lemma~\ref{Claim:1}}

Let $s \in \{1, \ldots, n\}$. We claim that it is enough to prove Lemma~\ref{Claim:1} for $t = s-1$.  Indeed, let $t \in \{0, \ldots, s-2\}$. Then, by the Alexandrov--Fenchel inequality and the symmetry of mixed volumes, we have the following
\[
\frac{V(\mathscr K_{s-t}, \mathscr L_{t}, \mathscr B_{n-s})}{V(\mathscr K_{s-t-1}, \mathscr L_{t+1}, \mathscr B_{n-s})} \le \frac{V(\mathscr K_{s-t-1}, \mathscr L_{t+1}, \mathscr B_{n-s})}{V(\mathscr K_{s-t-2}, \mathscr L_{t+2}, \mathscr B_{n-s})}. 
\]
If $\mathcal Q:=(\mathscr K_{s-t-2}, \mathscr L_{t}, \mathscr B_{n-s})$, then the last inequality becomes 
\[
\frac{V(K, K, \mathcal Q)}{V(K, L, \mathcal Q)} \le \frac{V(K, L, \mathcal Q)}{V(L,L, \mathcal Q)}. 
\]
Using the Alexandrov--Fenchel inequality $s-t-1 \ge 1$ times, we obtain the following sequence of inequalities:
\begin{equation*}
    \begin{aligned}
        \frac{V(\mathscr K_{s-t}, \mathscr L_{t}, \mathscr B_{n-s})}{V(\mathscr K_{s-t-1}, \mathscr L_{t+1}, \mathscr B_{n-s})} &\le \frac{V(\mathscr K_{s-t-1}, \mathscr L_{t+1}, \mathscr B_{n-s})}{V(\mathscr K_{s-t-2}, \mathscr L_{t+2}, \mathscr B_{n-s})} \le \frac{V(\mathscr K_{s-t-2}, \mathscr L_{t+2}, \mathscr B_{n-s})}{V(\mathscr K_{s-t-3}, \mathscr L_{t+3}, \mathscr B_{n-s})} \le \ldots \\
        &\ldots \le \frac{V(\mathscr K_1, \mathscr L_{s-1}, \mathscr B_{n-s})}{V(\mathscr L_{s}, \mathscr B_{n-s})}. 
    \end{aligned}
\end{equation*}
If the claim is true for $t = s-1$ then
\[
\frac{V(\mathscr K_1, \mathscr L_{s-1}, \mathscr B_{n-s})}{V(\mathscr L_{s}, \mathscr B_{n-s})} \le 1. 
\]
This yields the required inequality for any $t \in \{0, \ldots, s-2\}$ together with the equality case.

Now we prove the case $t=s-1$. Let 
\[
\cC:= (\mathscr L_{s-1}, \mathscr B_{n-s}).
\]
We need to show that
\begin{equation}
\label{Eq:Go}
\frac{V(\mathscr K_{1}, \mathscr L_{s-1}, \mathscr B_{n-s})}{V(\mathscr L_{s}, \mathscr B_{n-s})} = \frac{V(K, \mathscr L_{s-1}, \mathscr B_{n-s})}{V(L, \mathscr L_{s-1}, \mathscr B_{n-s})} = \frac{V(K, \mathcal C)}{V(L, \mathcal C)} \le 1,
\end{equation}
with equality if and only if $K$ is a $1$-convex lens.

Similar to the proof of Theorem~\ref{Thm:Main}, let $F_i$ be a facet of the $1$-convex polytope $K$ tangent to $B$ at $p_i$. Let $\mathcal N_{F_i} \subset  \mathbb S^{n-1}$ be the radial projection of $F_i$ onto the unit sphere $\mathbb S^{n-1}$ with respect to the origin. 

Assume that the $1$-convex lens $L$ is placed so that its center is at the origin and one of its facets, which we call $F$, touches the inscribed ball $B$ at the point $p_i$. We denote $\mathcal N_F \subset \mathbb S^{n-1}$ a similar radial projection of $F$ on $B$. We have the same inclusion~\eqref{Eq:Inclusion}.

The following lemma is an analogue of Lemma~\ref{Claim:M}. We postpone its proof to the end of this section. The main difference to the earlier proof is that, a priori, the measure $S(\sigma_i \cC)(u)$ might not have a smooth monotone density, and hence we cannot use the Fubini argument that we used in the proof of Lemma~\ref{Claim:M}. Instead, our argument will follow the approximation strategy similar to the one in \cite[Section 4]{DT}. 

    \begin{lemma}
    \label{Claim:M2}
        \begin{equation}
        \label{Eq:M}
        \frac{\int \limits_{\mathcal N_{F_i}} h_K(u) \, dS(\sigma_i \cC)(u)}{\int \limits_{\mathcal N_{F_i}} dS(\sigma_i \cC)(u)} \le \frac{\int\limits_{\mathcal N_F} h_L(u) \, dS(\sigma_i \cC)(u)}{\int \limits_{\mathcal N_F} dS(\sigma_i \cC)(u)} = \frac{V(L, \cC)}{V(\B, \cC)}.
    \end{equation}
    where $\sigma_i$ is the rotation about the origin that moves $L$ so that one of its faces touches the inscribed ball $B$ at $p_i$. Moreover, the equality holds if and only if $\mathcal N_F = \mathcal N_{F_i}$ (and thus, if and only if $F = F_i$). 
    \end{lemma}

By definition and our choice of $p_0$, each $\sigma_i$ is an element of $\mathscr{O}(K)$ that rotates $p_0$ to $p_i$. We put $\sigma_0=id$.

We know that $\sigma_i(K)=K$, and thus $h_{\sigma_i(K)} = h_K$ and $\mathcal N_{F_i} = \mathcal N_{\sigma_i(F_i)}$. Hence, by this invariance and \eqref{Eq:Isometry},
    \begin{equation}
        \label{Eq:B}
    \begin{aligned}
        \int \limits_{\mathcal N_{F_i}} h_K(u) \, dS(\cC)(u) = \int \limits_{\mathcal N_{\sigma_i(F_i)}} h_{\sigma_i(K)}(u) & \, dS(\sigma_i\cC)(u) = \int \limits_{\mathcal N_{F_i}} h_K(u) \, dS(\sigma_i \cC)(u),\\
        \int \limits_{\mathcal N_{F_i}} dS(\cC)(u) &=  \int \limits_{\mathcal N_{F_i}} dS(\sigma_i \cC)(u).
    \end{aligned}
    \end{equation}

Then, summing over all facets of $K$ and using \eqref{Eq:B} and Lemma~\ref{Claim:M2}, we obtain:
\begin{equation*}
    \begin{aligned}
        V(K, \cC) &= \frac{1}{n}\int \limits_{\mathbb S^{n-1}} h_K(u) \, dS(\cC)(u) = \frac{1}{n}\sum_i \int \limits_{\mathcal N_{F_i}} h_K(u) \, dS(\cC)(u) = \frac{1}{n}\sum_i \int \limits_{\mathcal N_{F_i}} h_K(u) \, dS(\sigma_i \cC)(u) \\
        &\le \frac{V(L, \cC)}{V(\B, \cC)} \cdot \frac{1}{n}\sum_i \int \limits_{\mathcal N_{F_i}} dS(\sigma_i \cC)(u) = \frac{V(L, \cC)}{V(\B, \cC)} \cdot \frac{1}{n}\sum_i \int \limits_{\mathcal N_{F_i}} dS(\cC)(u) \\
        &=\frac{V(L, \cC)}{V(\B, \cC)} \cdot \frac{1}{n}\int_{\mathbb S^{n-1}} dS(\cC)(u) = V(L, \cC), 
    \end{aligned}
\end{equation*}
which yields~\eqref{Eq:Go} and finishes the proof of Lemma~\ref{Claim:1}. 

\subsection{Proof of Lemma~\ref{Claim:M2}}

For brevity, let us denote $\mu := S(\sigma_i \mathcal C)$. Similar to Lemma~\ref{Claim:M}, it is enough to establish
\begin{equation}
    \label{Eq:M2}
    \frac{\int \limits_{\mathcal N_{F_i}} h_L(u) \, d\mu(u)}{\int \limits_{\mathcal N_{F_i}} d\mu(u)} \le \frac{\int\limits_{\mathcal N_F} h_L(u) \, d\mu(u)}{\int \limits_{\mathcal N_F} d\mu(u)} =: A_n.
\end{equation}
We need to show that
\begin{equation}
\label{Eq:Goal2}
\mathcal R[X] := \int\limits_{X} \left(h_L(u) - A_n\right) \,du \le 0 \quad \text{ for } \quad X = \mathcal N_{F_i} \subset\S^{n-1}.
\end{equation}
Note that $\mathcal R[\mathcal N_F]=0$. Furthermore, by construction, $\mathcal N_F$ is equal to the closed hemisphere in $\S^{n-1}$ centered at $p_i$, while $\mathcal N_{F_i}$ is a closed convex polytope in this hemisphere.   

Let us denote $g(u) := h_L(u) - A_n$. The function $g$ and the measure $\mu$ are both rotationally invariant with respect to the rotations in $\S^{n-1}$ about $p_i$. 

For $\theta \in [0,1]$, consider a pair of oriented hyperplanes passing through $op_i$ and making the (oriented) angle $2\pi \theta$. Denote by $H_\theta$ the wedge of angle $2\pi\theta$ between these hyperplanes. Finally, define $C_\theta := H_{\theta} \cap \mathcal N_F$. Pick $k \in \mathbb N_{>0}$ and let $C_{1/k}^j$, $j\in \{1, \ldots, k\}$, be a cyclically ordered set of $k$ spherical wedges with angle $2\pi/k$ that tile $\mathcal N_F$. We have
\[
\mathcal R[\mathcal N_F] = \int_{\mathcal N_{F}} g(u) d\mu(u)=0, \quad \bigcup_{j=1}^k C^j_{1/k} = \mathcal N_F.
\]
Hence, by rotational symmetry, 
\begin{equation}
\label{Eq:Slice}
\int_{C_{1/k}} g(u) d\mu(u) = 0, \quad \forall k \in \N_{>0}, \forall C_{1/k}.
\end{equation}

The function $g(u)$ is strictly monotonically increasing along any geodesic ray in $\mathcal N_F$ starting at $p_i$. The same is true for the function $h_L(u)$. Therefore, along any ray starting at $p_i$ there exists a unique point $u_0 \in \mathcal N_F$ where $g(u_0)=0$. This easily follows from the fact that $g(u)$ is rotationally invariant and $\mathcal R[\mathcal N_F] = 0$. All such points $u_0$ form an $(n-2)$-dimensional sphere $S^0 \subset \mathcal N_F$ centered at $p_i$ with the spherical radius $r_0$ satisfying $g(r_0) = 0$ (or equivalently,  $h_L(r_0) = A_n$). This sphere splits $\mathcal N_F$ into two regions: the open ball $D^-$ containing $p_i$, and a spherical annulus $D^+$, so that $\mathcal N_F = D^- \sqcup S_0 \sqcup D^+$. The function $g(u)$ is negative for every $u \in D^-$ and is positive for every $u\in D^+$. The same conclusion is true for $C_{1/k} \cap D^-$ and $C_{1/k} \cap D^+$. Therefore, for every geodesic ball $D_t \subset \S^{n-1}$ of spherical radius $t$ centered at $p_i$ we have that 
\begin{equation}
\label{Eq:Appriximation}
\int_{D_t \cap C_{1/k}} g(u) d \mu(u) \le 0, \quad \text{with equality if and only if} \quad t = \frac{\pi}{2}.
\end{equation}
Indeed, if $D_t \cap C_{1/k} \subset (D^- \sqcup S^0)\cap C_{1/k}$, then the conclusion is obviously true because $g$ is negative in $D^-$. If $D_t \cap C_{1/k} \supset (D^- \sqcup S^0) \cap C_{1/k}$ and $\int_{D_t \cap C_{1/k}} g(u) d \mu(u) \ge 0$, then 
\[
\int_{C_{1/k}} g(u) d \mu(u) = \int_{D_t \cap C_{1/k}} g(u) d\mu(u) + \int_{(D^+ \setminus D_t) \cap C_{1/k}} g(u) d\mu(u) > 0
\]
because $g$ is positive in $D^+$; this is a contradiction to \eqref{Eq:Slice}. Hence \eqref{Eq:Appriximation} is true.

Construct a sequence $C^j_{1/k}$ of wedges, $j \in \{1,\ldots, k\}$, such that $C^j_{1/k}$'s have disjoint interiors, and a sequence $D_{t_{k,j}}$ of balls such that for every $k \in \N$,  
\[
\mathcal N_F = \bigcup_{j=1}^k C^j_{1/k}, \quad \quad 
\mathcal N_{F_i} \subset \bigcup_{j=1}^k \big(D_{t_{k,j}} \cap C^j_{1/k}\big),  
\]
and such that
\[
\bigcup_{j=1}^k \big(D_{t_{k,j}} \cap C^j_{1/k}\big) \underset{k\to \infty}{\longrightarrow} 
\mathcal N_{F_i}
\]
in the Hausdorff metric on $\S^{n-1}$. More specifically, for each $C^j_{1/k}$, we define $t_{k,j}$ to be the smallest radius such that the ball $D_{t_{k,j}}$ contains $\mathcal N_{F_i} \cap C^j_{1/k}$. This approximation is similar to the one in \cite{DT}.

Therefore, since the function $g$ is bounded on $\mathcal N_F$, we obtain 
\[
\int \limits_{\mathcal N_{F_i}} g(u)\,d\mu(u)=\lim_{k \to \infty} \sum_{j=1}^k \int \limits_{D_{t_{k,j}} \cap C^j_{1/k}} g(u)\, d\mu(u) \le^{\eqref{Eq:Appriximation}} 0,
\]
which is the desired inequality \eqref{Eq:Goal2}.

Let us analyze the equality case. If $\mathcal N_{F_i} \subset D^- \sqcup S^0$, then $\int_{\mathcal N_{F_i}} g(u) d\mu(u) < 0$. Thus, we can assume that $\mathcal N_{F_i}$ intersects $D^+$. If $\mathcal N_F\neq \mathcal N_{F_i}$, then there exists a sufficiently large $k_0$ and $j_0 \in \{1, \ldots, k_0\}$, such that 
\begin{equation}
\label{Eq:Inc}
\mathcal N_F \cap C^{j_0}_{1/k_0} \supsetneq D_{t_{k_0,j_0}} \cap C^{j_0}_{1/k_0} \supseteq\mathcal N_{F_i} \cap C^{j_0}_{1/{k_0}} \supsetneq D^- \cap C^{j_0}_{1/k_0}.
\end{equation}
For $k \ge k_0$, we can split the union $\mathcal U_k:=\bigcup_{j=1}^k \big(D_{t_{k,j}} \cap C^j_{1/k}\big)$ into two subsets $\mathcal X_k$ and $\mathcal Y_k$ with disjoint interiors such that 
\[
\mathcal X_k := \overline{C^{j_0}_{1/k_0} \cap \mathcal U_k}, \quad \mathcal Y_k := \overline{\mathcal U_k \setminus \mathcal X_k}.
\]
By passing to a subsequence if necessary, we can assume that for $k \ge k_0$ if $C^{j}_{1/k} \cap C^{j_0}_{1/k_0} \neq \emptyset$, then $C^{j}_{1/k} \subseteq C^{j_0}_{1/k_0}$ (for example, we can choose $k = 2^s \cdot k_0$ for $s=1,2,\ldots$). In this way, for $k \ge k_0$, we have $\int_{\mathcal Y_k} g(u) d\mu(u) \le 0$ and
\begin{equation*}
    \begin{aligned}
        \sum_{j=1}^k \int \limits_{D_{t_{k,j}} \cap C^j_{1/k}} g(u)\, d\mu(u) &= \int \limits_{\mathcal X_k} g(u)\, d\mu(u) + \int \limits_{\mathcal Y_k} g(u)\, d\mu(u) \\
        &\le^{\eqref{Eq:Inc}} \int \limits_{D_{t_{k_0,j_0}} \cap C^{j_0}_{1/k_0}} g(u) \, d\mu(u) =: E <^{\eqref{Eq:Appriximation}\eqref{Eq:Inc}} 0        
    \end{aligned}
\end{equation*}

Therefore, $\int \limits_{\mathcal N_{F_i}} g(u)\,d\mu(u)=\lim_{k \to \infty} \sum_{j=1}^k \int \limits_{D_{t_{k,j}} \cap C^j_{1/k}} g(u)\, d\mu(u) \le E < 0$, thus we have the strict inequality. Thus we obtained that the only way we can achieve equality is when $\mathcal N_F = \mathcal N_{F_i}$.  This finishes the proof of Lemma~\ref{Claim:1}.


\begin{thebibliography}{FIMP2}

\bibitem[AACF]{AACF}
S. Artstein-Avidan, A. Chor, D. Florentin, \emph{A full classification of the isometries of the class of ball-bodies}, 2025, \url{https://arxiv.org/abs/2503.02613}.

\bibitem[AAF]{AAF} S. Artstein-Avidan, D. Florentin. \emph{An in-depth study of ball-bodies}, 2025, \url{https://arxiv.org/abs/2505.09200}.

\bibitem[Be]{Bezdek2008}
 K. Bezdek, \emph{From the Kneser–Poulsen conjecture to ball-polyhedra}, European Journal of Combinatorics, {\bf29}(8) (2008), 1820--1830.

\bibitem[Be2]{Bezdek2012} 
K. Bezdek,  \emph{Illuminating spindle convex bodies and minimizing the volume of spherical sets of constant width}, Discret. Comput. Geom. \textbf{47}(2) (2012), 275--287.

\bibitem[Be3]{BezdekConj}
K. Bezdek, \emph{Volumetric bounds for intersections of congruent balls in Euclidean spaces}. Aequationes Math. \textbf{95} (2021), no. 4, 653--665.

\bibitem[BLN]{BLN}
K. Bezdek, Z. L\'angi, M. Nasz\'odi, \emph{Selected topics from the theory of intersections of balls}, 2024, \url{https://arxiv.org/abs/2411.10302}.

\bibitem[BLNP]{BLNP}
K. Bezdek, Z. L\'angi, M. Nasz\'odi, P. Papez, \emph{Ball-polyhedra}. Discrete Comput. Geom. \textbf{38} (2007), no. 2, 201--230.

\bibitem[Bla]{Bla56} W. Blaschke,
\emph{Kreis und Kugel}. de Gruyter, Berlin, 1956.

\bibitem[BD1]{BorDr13} A. Borisenko, K. Drach, \emph{Closeness to spheres of hypersurfaces with normal curvature bounded below}. Sb. Math. 204:11 (2013), 1565--1583.

\bibitem[BD2]{BD} A. Borisenko, K. Drach, \emph{Isoperimetric inequality for curves with curvature bounded below}. Math. Notes \textbf{95} (2014), no. 5, 590--598.

\bibitem[BD3]{BorDr15_1} A. Borisenko, K. Drach, \emph{Extreme properties of curves with bounded curvature on a sphere}. J. Dyn. Control Syst. \textbf{21} (2015), no. 3, 311--327.

\bibitem[BM]{BM} A. Borisenko, V. Miquel, \emph{Total curvatures of convex hypersurfaces in hyperbolic space}. Illinois J. Math. \textbf{43} (1999), no. 1, 61--78.

\bibitem[BH]{BH}  K.~J. B\"or\"oczky, D. Hug, \emph{A reverse Minkowski-type inequality}. Proc. Amer. Math. Soc. {\bf 48} (2020), 4907--4922.

\bibitem[Dr1]{DrSpindle}
K. Drach, \emph{Some sharp estimates for convex hypersurfaces of pinched normal curvature}. J. Math. Phys. Anal. Geom. \textbf{12} (2015), no. 2, 111--122.

\bibitem[Dr2]{Dr}
K. Drach, \emph{Inradius estimates for convex domains in 2-dimensional Alexandrov spaces}. Anal. Geom. Metr. Spaces \textbf{6} (2018), no. 1, 165--173. 

\bibitem[Dr3]{DrBla}
K. Drach, \emph{The Blaschke rolling theorem in Riemannian manifolds of bounded curvature}, 2024, \url{https://arxiv.org/abs/2404.02739}.

\bibitem[DT]{DT}
K. Drach, K. Tatarko, \emph{Reverse isoperimetric problems under curvature constraints}, 2023, \url{https://arxiv.org/abs/2303.02294}.

\bibitem[FV]{FV}F. Fodor,  V. V\'igh, \emph{Disc-polygonal approximations of planar spindle convex sets}, Acta Sci. Math. (Szeged) {\bf78} (2012), no. 1-2, 331--350.

\bibitem[FKV]{FKV}  F. Fodor, A. Kurusa, V. V\'igh, \emph{Inequalities for hyperconvex sets}. Adv. Geom. {\bf16} (2016), no. 3, 337--348.

\bibitem[Li]{Li1977} J. Linhart, \emph{Kantenl\"angensumme, mittlere Breite und Umkugelradius konvexer Polytope}, Arch. Math. {\bf29}(5)  (1977), 558--560.

\bibitem[Mi]{MilInradius}
A. Milka, \emph{Estimates of the sizes of curves with bounded curvature}. Ukrainian Geom. Sb. \textbf{21} (1978), 88--91. 

\bibitem[NV]{NV}
K. Nagy, V. V\'igh, \emph{Best and random approximations with generalized disc-polygons}, Discrete Comput. Geom. {\bf72} (2024), no. 1, 357--378. 

\bibitem[Sch]{Sch} R. Schneider, \emph{Convex bodies: The Brunn-Minkowski theory}. Second expanded edition, Cambridge Press, 2014.

\bibitem[SWY]{SWY} C. Sch\"utt, E. Werner, D. Yalikun, \emph{Floating bodies for ball-convex bodies}, 2025. \url{https://arxiv.org/abs/2504.15488}.

\end{thebibliography}
\end{document}